\documentclass[11pt, leqno]{amsart}
\usepackage{anysize}
\usepackage[activeacute,english]{babel}%
\usepackage[ansinew]{inputenc}
\usepackage{graphicx}
\usepackage{amsmath}
\usepackage{amsfonts}
\usepackage{amssymb}
\usepackage{amsxtra}
\usepackage{amstext}
\usepackage{amscd}
\usepackage{latexsym}
\usepackage{stmaryrd}
\usepackage{mathrsfs}
\usepackage{multicol}
\usepackage{fancyhdr}
\usepackage[mathscr]{euscript}

\usepackage[hidelinks]{hyperref}

\makeatletter
\newtheorem{teo}{Theorem}[section]

\newtheorem{defi}[teo]{Definition}

\newtheorem{rem}{Remark}

\begin{document}

\title[A note on one-dimensional symmetry for Hamilton-Jacobi equations]{A note on one-dimensional symmetry for Hamilton-Jacobi equations with extremal Pucci operators and application to Bernstein type estimate}
\author[Rodrigo Fuentes]{R. Fuentes$^{*}$}
\author[Alexander Quaas]{A. Quaas$^{\dag}$}

\address[$\dag$]{Departamento de  Matem\'atica, Universidad T\'ecnica Federico Santa Mar\'ia, Casilla 110-V, Valpara\'iso, Chile.}
\email[$1$]{rfuentes.montecinos@gmail.com}
\email[$2$]{alexander.quaas@usm.cl}
\date{\today}
\maketitle

\begin{abstract}
{{We prove a Liouville-type theorem that is one-dimensional symmetry and classification results for non-negative $L^q$-viscosity solutions of the equation
\begin{equation*}
-\mathcal{M}_{\lambda, \Lambda}^{\pm}(D^2u)\pm |Du|^p=0, x\in \mathbb{R}_+^n,
\end{equation*}
with boundary condition $u(\tilde{x},0)=M\geq 0, \tilde{x}\in \mathbb{R}^{n-1}$, where $\mathcal{M}_{\lambda, \Lambda}^{\pm}$ are the Pucci's operators with parameters $\lambda, \Lambda \in \mathbb{R}_+$ $0<\lambda\leq \Lambda$ and $p>1$. 
The results are an extension of the results by Porreta and Ver\'on in \cite{porreta} for the case  $p\in (1,2]$ and by o Filippucci, Pucci and Souplet in \cite{fps} for the case $p>2$, both for the Laplacian case (i.e. $\lambda=\Lambda=1$).
As an application in the case $p>2$, we prove a sharp Bernstein estimation for $L^q$-viscosity solutions of the fully nonlinear equation 
\begin{equation*}
-\mathcal{M}_{\lambda, \Lambda}^{\pm}(D^2u)= |Du|^p+f(x), \quad x\in \Omega, \label{ecuacion1}
\end{equation*}
with boundary condition $u=0$ on $\partial \Omega$, where $\Omega \subset \mathbb{R}^n$. }}
\end{abstract}

\section{Introduction}
\noindent In this paper we consider, for $n\geq 2$, the second order elliptic equation
\begin{equation} \label{ecuacioneliptica}
\begin{cases}
-\mathcal{M}_{\lambda, \Lambda}^{\pm}(D^2u)\pm |Du|^p=0, &x\in \mathbb{R}_+^n,\\
u(\tilde{x},0)=M,& \tilde{x}\in \mathbb{R}^{n-1},\\
u\geq 0, & \mbox{ in } \mathbb{R}_+^n
\end{cases}
\end{equation}
where $\mathbb{R}_+^n=\left\lbrace (\tilde{x},x_n)=(x_1,\cdots ,x_n)\in \mathbb{R}^n,x_n>0\right\rbrace $, $p>1$, $M\geq 0$ and $\mathcal{M}_{\lambda,\Lambda}^{\pm}$ are the extremal Pucci's operators with parameters $0\leq \lambda \leq \Lambda$ defined by
\begin{equation}\label{definicionpucci}
\begin{array}{cc}
\mathcal{M}_{\lambda,\Lambda}^+ (N)&=\Lambda \displaystyle\sum_{e_i >0}e_i+\lambda \displaystyle\sum_{e_i<0} e_i,\\
\mathcal{M}_{\lambda,\Lambda}^- (N)&=\lambda \displaystyle\sum_{e_i >0}e_i+\Lambda \displaystyle\sum_{e_i<0} e_i,
\end{array}
\end{equation}
with $N$ any symmetric $n\times n$ matrix, and $e_i$ the eigenvalues of $N$. 
One first aim of this paper is to establish one-dimensional symmetry for any $L^q$-viscosity solution of equation (\ref{ecuacioneliptica}).
We split the results in two cases depending on $p$. For $p\in (1,2]$ we have  

\begin{teo} \label{teoremapmenora2}
Let $p\in (1,2]$ and $u\in C(\overline{\mathbb{R}_+^n})$ an $L^q$-viscosity solution of
\begin{equation} \label{ecuacioneliptica1}
\begin{cases}
-\mathcal{M}_{\lambda, \Lambda}^{\pm}(D^2u)+ |Du|^p=0, &x\in \mathbb{R}_+^n,\\
u(\tilde{x},0)=M,& \tilde{x}\in \mathbb{R}^{n-1},\\
u\geq 0, & \mbox{ in } \mathbb{R}_+^n,
\end{cases}
\end{equation}
with $q\geq n$ and $M\geq 0$. Then $u$ depends only on the variable $x_n$.
\end{teo}

and for $p>2$ we have
\begin{teo} \label{teoremapmayora2}
Let $p>2$ and $u\in C(\overline{\mathbb{R}_+^n})$ an $L^q$-viscosity solution of
\begin{equation} \label{ecuacioneliptica2}
\begin{cases}
-\mathcal{M}_{\lambda, \Lambda}^{\pm}(D^2u)=|Du|^p, &x\in \mathbb{R}_+^n,\\
u(\tilde{x},0)=0,& \tilde{x}\in \mathbb{R}^{n-1},\\
u\geq 0, & \mbox{ in } \mathbb{R}_+^n,
\end{cases}
\end{equation}
with $q\geq n$. Then $u$ depends only on the variable $x_n$.
\end{teo}

From these theorems, classification results can be established by solving the ODE, see the beginning of section 3. 
Notice that the only work we know of one-dimensional symmetry or rigidity results in the half space for fully non-linear operator is in \cite{leoni}, for explosive boundary condition. Other symmetry type results in bounded domain can be found in \cite{boyan}.

In the case of the Laplacian rigidity results in the half-space is widely studied for $-\Delta u=f(u)$ in different setting see for example 
\cite{bere}, \cite{bere1} and \cite{jorge}, other references can be founded in these works.

 In the context of our problem and the Laplacian, Porreta and Ver\'on in \cite{porreta} studied case $p\in (1,2]$ and Filippucci, Pucci and Souplet 
\cite{fps}, that are connected with our setting.

One useful tool in establishing the above theorems is a Bernstein type estimate by Birindelli, Demengel and Leoni in \cite{bernstein}, specifically Proposition 2.3. with $\alpha=0$, which is an extension of the result from Capuzzo Dolcetta, Leoni and Porreta in \cite{clp}. 

This Bernstein-type result are for second order elliptic equations of the form
\begin{equation}
-\mathcal{M}_{\lambda, \Lambda}^{\pm}(D^2u) +g(Du)=f(x), \quad x\in \Omega, \label{ecuacion1}
\end{equation}
where $\Omega$ is an open bounded subset of $\mathbb{R}^n$.\\

Using a regularity argument and the Bernstein estimate, we prove first that any $L^q$-viscosity solution of problem (\ref{ecuacioneliptica}) has to be a classic solution. 

Now we come back to the proof of our one-dimensional symmetry results. For $p\in (1,2]$ the proof is based in the fact that the solution of the problem has a finite limit when de variable $x_n$ goes to infinity. In \cite{porreta} this fact is guaranteed by the Bernstein estimation from \cite{lions}, and in this paper the result is obtained thanks to the Theorem \ref{teoremabernstein} in section 4. For $p>2$, we use a moving planes argument and the maximum principle to prove that any $L^q$-viscosity solution of (\ref{ecuacioneliptica}) is one dimensional. The moving planes technique was used in \cite{fps} to prove the same result for a classic solution to the equation
$$\begin{cases}
-\Delta u=|Du|^p, &x\in \mathbb{R}_+^n,\\
u(\tilde{x},0)=0,& \tilde{x}\in \mathbb{R}^{n-1},
\end{cases}$$
using the Bernstein estimation in \cite{lions}, and a regularity argument.

\bigskip

\noindent We present here for completeness in section 4 a Bernstein estimate that we will use through the paper. This result is a simplification of the work done in \cite{bernstein}, and we present a shorter proof, that is based on the same ideas of \cite{clp} as the proof in \cite{bernstein}, that works also for the case $\lambda=0$.

  The Bernstein technique was implemented for the first time by Bernstein himself in \cite{bernstein0}, and it has been replicated for equations with different elliptical operators. For example, in 1985, Lions in \cite{lions} proved for the equation
\begin{equation} \label{ecuacionejemplo2}
-\Delta u+\alpha|Du|^p=f(x),\quad \mbox{ in } \Omega,
\end{equation} 
with $\alpha>0$, the estimation $$|Du|\leq C \left\lbrace d(x,\partial \Omega)^{-\frac{1}{p-1}}+C_1\right\rbrace ,$$ where $C_1$ depends on $f$, for $p>1$. Lions used that $f\in C^1(\Omega)$, and the proof was based on the traditional Bernstein technique, using the equation for the function $|Du|^2$. Many other Bernstein type results have been proof for other elliptic operators. In \cite{barles}, Barles extended this idea to a weak Bernstein method to prove global Lipschitz continuity, using only continuous viscosity solutions and bases his work not in the fact that derive the equation, but only in properties of the operator $F$. More recently in \cite{barles2}, Barles use this weak Bernstein method to prove local Lipschitz continuity for a continuous solution of 
\begin{equation*} \label{ecuacionejemplo}
F(x,u,Du,D^2u)=0 \quad \mbox{ in } \Omega,
\end{equation*} 
using locally Lipschitz continuity and ellipticity, with other properties of $F$. Others regularity results for viscosity solutions were proved by Kawohl and Kutev in \cite{kawohl}, and by Armstrong and Tran in \cite{armstrong} using different  assumptions. In \cite{kawohl} they used a boundary Lipschitz condition (and more recently in \cite{kawohl2} using $C^2$-boundary condition) to prove Lipschitz regularity in all the domain $\overline{\Omega}$, and in \cite{armstrong} they require that the diffusion matrix $A$ from the equation 
$$-tr(A(x)D^2u)+H(x,Du)=0, \quad \mbox{ in } B_2$$
to be regular enough in view to prove that any continuous solution is Lipschitz continuous in $B_1$. H$\ddot{\mbox{o}}$lder regularity results for extensions of equation (\ref{ecuacionejemplo2}) have recently been published in \cite{birindelli} and \cite{dasilva}. Da Silva and Nornberg use in \cite{dasilva} $L^q$-viscosity solutions for elliptic equations with Hamiltonian with coefficients not necessarily regular, as in our equation (\ref{ecuacion1}), to prove local H$\ddot{\mbox{o}}$lder continuity.\\

As an application to their Liouville result, Filipucci, Pucci and Souplet proved an improvement of the Bernstein estimation in \cite{lions}, now with the boundary value, for the inhomogeneous Dirichlet problem (specifically Theorem 1.6 in \cite{fps}). The result establishes a more precise constant in the Bernstein estimate, which is also valid for  equation of the form
\begin{equation} \label{ecuacionelipticanohomogenea}
\begin{cases}
-\mathcal{M}_{\lambda, \Lambda}^{\pm}(D^2u)=|Du|^p+f(x), \quad \mbox{in}\,\, \Omega,\\
u=0, \quad \mbox{on} \,\, \partial \Omega,
\end{cases}
\end{equation}
this by using our Theorem \ref{teoremapmayora2} and similar ideas as in \cite{fps}.

Specifically, we have
\begin{teo} \label{teoremabernsteinmejorado}
Let $p>2$ and $f$ a Lipschitz function with $\|f\|_{W^{1,\infty}}\leq M$. If $u\in C(\overline{\Omega})$ is an $L^q$-viscosity solution of (\ref{ecuacionelipticanohomogenea}), then for every $\epsilon >0$ there exists a constant $C=C(\epsilon, M)>0$ such that
$$|Du(x)|\leq (1+\epsilon)\left( \dfrac{\lambda}{p-1} \right)^{\frac{1}{p-1}}d_{\partial \Omega}(x)^{-\frac{1}{p-1}}+C \quad \mbox{in} \,\, \Omega,$$
if the operator in (\ref{ecuacionelipticanohomogenea}) is $\mathcal{M}^+$. If the operator is $\mathcal{M}^-$, the constant must be $\Lambda$.
\end{teo}

The paper is organized as follows: In Section 2 we give some preliminaries; in Section 3 we prove our main results; in section 4 we give and prove our first Bernstein type estimate.
\section{Preliminaries}

\noindent Next we are going to give the different definitions of solutions that we will use throughout the paper. For the definitions, with $\Omega \subset \mathbb{R}^n$, we consider a function $F:\Omega \times \mathbb{R}\times \mathbb{R}^n \times \mathcal{S}(n) \rightarrow \mathbb{R}$, where $\mathcal{S}(n)$ denotes the set of $n\times n$ real symmetric matrices, and the equation
\begin{equation} \label{ecuaciondef}
F(x,u(x),Du,D^2u)=f(x), \quad \mbox{ in } \Omega.
\end{equation}
\begin{defi}
Let $n\leq q$ and $f\in L_{loc}^q(\Omega)$. A function $u \in C(\Omega)$ is an $L^q$-viscosity subsolution (supersolution) of (\ref{ecuaciondef}) if for all $\varphi \in W_{loc}^{2,q}(\Omega)$, and $\hat{x} \in \Omega$ at which $u-\varphi$ has a local maximum (respectively minimum) we have
\begin{equation*}
\begin{cases}
\mathrm{ess}\, \displaystyle\liminf_{x\rightarrow \hat{x}} (F(x,u(x),D\varphi(x),D^2\varphi(x))-f(x))\leq 0\\
(\mathrm{ess}\, \displaystyle\limsup_{x\rightarrow \hat{x}} (F(x,u(x),D\varphi(x),D^2\varphi(x))-f(x))\geq 0).
\end{cases}
\end{equation*}
In addition, $u$ is an $L^q$-viscosity solution of (\ref{ecuaciondef}) if it is both an $L^q$-viscosity subsolution and an $L^q$-viscosity supersolution.
\end{defi}
\noindent This notion of solution was studied by Caffarelli, Crandall, Kocan and Swiech in \cite{measurable}, and rightly related with the notion of $L^q$-strong solutions. The use of $L^q-$viscosity solutions falls in that not differentiability is required for the test functions, either the continuity of $F$ and $f$ at $x$, therefore, pointwise relations are not required. If the functions $F$ and $f$ are continuous, and the test function $\varphi$ is twice continuously differentiable we employ the notion of classical viscosity solution that we define next.
\begin{defi}
A function $u\in C(\Omega)$ is a viscosity subsolution (supersolution) of (\ref{ecuaciondef}) if, for any function $\varphi\in C(\Omega)\cap C^2(\Omega)$, if $\hat{x}$ is a global maximum (minimum) point of $u-\varphi$, then
\begin{equation*}
\begin{cases}
F(\hat{x},u(\hat{x}),D\varphi(\hat{x}),D^2\varphi(\hat{x}))\leq f(\hat{x})\\
(F(\hat{x},u(\hat{x}),D\varphi(\hat{x}),D^2\varphi(\hat{x}))\geq f(\hat{x})).
\end{cases}
\end{equation*}
In addition, $u$ is a viscosity solution of (\ref{ecuaciondef}) if it is both a viscosity subsolution and a viscosity supersolution.
\end{defi}
\noindent The use of these types of solutions allows establish results for uniqueness, existence and stability for a wide class of equations. The notion of viscosity solutions was introduced by Crandall and Lions in \cite{viscosity} for first-order equations of Hamilton-Jacobi type. A second-order extension was introduced by Lions in \cite{lions2} for Hamilton-Jacobi-Bellman equations, where he showed uniqueness results for viscosity solutions of these type of problems. \\

\noindent We use that the Pucci operators are extremal, and can be written by
\begin{equation}
\begin{array}{cc} \label{definicionsup}
\mathcal{M}_{\lambda,\Lambda}^+ (N)&=\displaystyle\sup_{A\in \mathcal{A}_{\lambda,\Lambda}} tr(AN),\\
\mathcal{M}_{\lambda,\Lambda}^- (N)&=\displaystyle\inf_{A\in \mathcal{A}_{\lambda, \Lambda}} tr(AN),
\end{array} 
\end{equation}
where $\mathcal{A}_{\lambda,\Lambda}$ denotes the set of all symmetric matrix whose eigenvalues lie in the interval $[\lambda,\Lambda]$. With these definitions, it is easy to see that the Pucci operators $\mathcal{M}_{\lambda, \Lambda}^{\pm}$ have the following properties: for any $\alpha\geq 0$ and for any $n\times n$ symmetric matrices $M,N$ we have
\begin{equation} \label{elliptic}
\begin{cases}
\mathcal{M}_{\lambda, \Lambda}^{\pm}(\alpha M)=\alpha \mathcal{M}_{\lambda,\Lambda}^{\pm}(M),\\
\mathcal{M}_{\lambda,\Lambda}^- (M) \leq \mathcal{M}_{\lambda,\Lambda}^{\pm} (M+N)-\mathcal{M}_{\lambda,\Lambda}^{\pm}(N) \leq \mathcal{M}_{\lambda,\Lambda}^+ (M).
\end{cases}
\end{equation}
(For these and other properties of Pucci operators see \cite{pucci}). 

\section{Proof of main results and classification solution}

In this section we prove our main results. 
We first start this section with the classification results that can be deduce from our main theorems.

 As a consequence of Theorem \ref{teoremapmenora2} we have that any solution of (\ref{ecuacioneliptica1}) is given by
\begin{enumerate}
\item[i)] For $p\in (1,2)$, $u\equiv 0$ or there exists a constant $l\in [0,M)$ such that
$$u(\tilde{x},y)=l+\displaystyle\int_y^{\infty} \left[ \left( \dfrac{p-1}{\Lambda }\right) (t+C_{M,l}) \right]^{-\frac{1}{p-1}}dt,$$
where $C_{M,l}$ is determined by the relation
$$\displaystyle\int_0^{\infty} \left[ \left( \dfrac{p-1}{\Lambda }\right)( t+C_{M,l})\right]^{-\frac{1}{p-1}}dt=M-l,$$
in the case of $\mathcal{M}_{\lambda, \Lambda}^+$. If the operator is $\mathcal{M}_{\lambda, \Lambda}^-$, the constant $\Lambda$ in the solution must be replaced by $\lambda$.
\item[ii)] For $p=2$ necessarily $u \equiv M$. 
\end{enumerate}

\noindent  As a consequence of Theorem \ref{teoremapmayora2} , in this case we have that any solution of (\ref{ecuacioneliptica2}) is given by $u\equiv 0$ or $u=u_{\hat{c}}(x_n)=u_0(x_n+\hat{c})-u_0(\hat{c})$ for some $\hat{c}\geq 0$, with
\begin{equation} \label{solucionfundamentalpmayora2}
u_0(t)=c_{p,\lambda,\Lambda}t^{\frac{p-2}{p-1}}, \quad t\geq 0,
\end{equation}
where $c_{p,\lambda,\Lambda}=\left( \dfrac{\lambda}{p-1} \right)^{\frac{1}{p-1}} \left( \dfrac{p-1}{p-2} \right)$ if the operator is $\mathcal{M}_{\lambda, \Lambda}^+$. In the case of $\mathcal{M}_{\lambda, \Lambda}^-$, the constant $\lambda$ in $c_{p,\lambda,\Lambda}$ must be replaced by $\Lambda$. \\

\noindent We note that in all of cases, the functions are classic solutions, this will be the first thing we will show next.

\bigskip

\noindent \textit{Proof of Theorem \ref{teoremapmenora2}:} First, we will show that $u$ is a classic solution of (\ref{ecuacioneliptica1}). By Theorem \ref{teoremabernstein} below, with $g(Du)=|Du|^p$, we have 
\begin{equation} \label{bernstein}
|Du(\tilde{x},y)|\leq C(n,p)y^{-\frac{1}{p-1}}, \quad a.e.\,\,(\tilde{x},y)\in \mathbb{R}^{n-1}\times (0,\infty).
\end{equation}
Considering $f(\tilde{x},y)=|Du(\tilde{x},y)|^p$ in Corollary 3.10. in \cite{measurable}, recalling that $\tilde{B}=B((\tilde{x},y),\frac{y}{2})$ satisfies the uniform exterior cone condition and $u$ is continuous in $\overline{\mathbb{R}_+^n}$, by (\ref{bernstein}) we have that $f \in L^q(\tilde{B})$ for all $n\leq q\leq \infty$, therefore $u\in W_{loc}^{2,q}(\mathbb{R}_+^n)$ (using by Lemma 2.5 in \cite{measurable} the fact that the unique $L^q$-strong solution of (\ref{ecuacioneliptica1}) coincides with our $L^q$-viscosity solution). This implies that $Du\in W_{loc}^{1,q}(\mathbb{R}_+^n)$, and by Morrey's inequality (Theorem 4 in 5.6.2, \cite{evans}), we have that $Du \in C_{loc}^\alpha(\mathbb{R}_+^n)$, for some $\alpha \in (0,1)$. Therefore, by Theorem 1.2 in \cite{cabre}, the function $u$ is in $C_{loc}^{2,\alpha}(\mathbb{R}_+^n)$, and this implies that $u$ is a classic solution of (\ref{ecuacioneliptica1}).\\\\
\noindent Now we will prove the symmetry result.\\\\
First we take $p \in (1,2)$, and we follow the ideas in the proof in \cite{porreta}. By Theorem \ref{teoremabernstein} and recalling that $u$ is a classic solution of (\ref{ecuacioneliptica1}) we have that
\begin{equation} \label{bernsteintodopunto}
|Du(\tilde{x},y)|\leq C(n,p)y^{-\frac{1}{p-1}}, \quad \forall \,\,(\tilde{x},y)\in \mathbb{R}^{n-1}\times (0,\infty).
\end{equation}
and therefore
\begin{equation} \label{cotaintegral}
|u(\tilde{x},\eta)-u(\tilde{x},y)|\leq C \displaystyle\int_{\eta}^y t^{-\frac{1}{p-1}}dt.
\end{equation}
Since $p\in (1,2)$, then $\frac{1}{p-1}>1$ and therefore $u(\tilde{x},y)$ has a finite limit as $y \rightarrow \infty$. By (\ref{bernsteintodopunto}) and using the mean value theorem, the limit does not depend on $\tilde{x}$, and then we can define $l:=\displaystyle\lim_{y\rightarrow \infty} u(\tilde{x},y)$, and by (\ref{cotaintegral}) we have that
\begin{equation} \label{cotasu}
l-Cy^{-\frac{2-p}{p-1}}\leq u(\tilde{x},y)\leq l+Cy^{-\frac{2-p}{p-1}}, \quad \forall (\tilde{x},y) \in \mathbb{R}_+^n.
\end{equation}
Now we want to prove that $u(\tilde{x},y)=v_l(y)$, where $v_l$ is the unique solution of the one-dimentional problem
\begin{equation} \label{onedimentional}
\begin{cases}
\Lambda v_l''&=|v'_l|^p \mbox{ in } (0,\infty),\\
v_l(0)&=M,\\
\displaystyle\lim_{y\rightarrow \infty} v_l(y)&=l.
\end{cases}
\end{equation}
We use here the constant $\Lambda$. If the operator is $\mathcal{M}_{\lambda, \Lambda}^-$, this constant is replaced by $\lambda$. To prove that $u\leq v_l$, for $t\in (0,1)$ and $C_R\in \mathbb{R}$, we consider the problem
\begin{equation} \label{ecuacionpsi}
\begin{cases}
-\mathcal{M}_{\Lambda,\lambda}^{+}(D^2\psi_{t,R})+\sqrt{1-t^2}|D\psi_{t,R}|^p+C_R=0 \quad \mbox{ in } B_R^{n-1}\subset \mathbb{R}^{n-1},\\
\psi_{t,R}(0)=0, \quad \displaystyle\lim_{|\tilde{x}| \rightarrow R^-} \psi_{t,R}(\tilde{x})=\infty.
\end{cases}
\end{equation}
We know by \cite{demengel} that there exists a unique constant $C_R$ such that the problem (\ref{ecuacionpsi}) has a solution $\psi_{t,R}$, and by \cite{birindelli}, this solution is also unique. Also by \cite{demengel} we know that $C_R>0$ and that $\psi_{t,R}$ achieves its minimum in zero, and writing $C_1$ and $\psi_{t,1}$ the solutions of (\ref{ecuacionpsi}) in $B_1^{n-1}$, we have
\begin{equation} \label{condicionr}
C_R=R^{-\frac{p}{p-1}}C_1, \quad \psi_{t,R}=R^{-\frac{2-p}{p-1}}\psi_{t,1}\left( \dfrac{|\tilde{x}|}{R}\right).
\end{equation}
Consider now, for $L\in \mathbb{R}_+$, the function $\varphi_{t,L,R}$ solution of the problem
\begin{equation} \label{ecuacionvarphi}
\begin{cases}
-\Lambda \varphi''+t|\varphi'|^p=\frac{\sqrt{1-t^2}}{t}C_R, \quad \mbox{ in } (0,L),\\
\varphi(0)=\frac{M}{t},\quad \varphi(L)=\frac{1}{t}\left( l+Cl^{-\frac{2-p}{p-1}} \right).
\end{cases}
\end{equation}
Define now the function $\overline{u}(\tilde{x},y)=t\varphi_{t,L,R}(y)+\sqrt{1-t^2}\psi_{t,R}(\tilde{x})$, we want to see that $\overline{u}$ is a supersolution of (\ref{ecuacioneliptica1}) in the cylinder $B_R^{n-1}\times (0,L)$. Indeed, using equations (\ref{ecuacionpsi}) and (\ref{ecuacionvarphi}) we have
\begin{align*}
-\mathcal{M}_{\Lambda,\lambda}^{+}(D^2\overline{u})+|D\overline{u}|^p&=t\left( \dfrac{\sqrt{1-t^2}}{t}C_R-t|\varphi_{t,L,R}'|^p \right)+\sqrt{1-t^2}\left(-\sqrt{1-t^2}|D\psi_{t,R}|^p-C_R\right)\\
&+\left( t^2|\varphi_{t,L,R}'|^2+(1-t^2)|D\psi_{t,R}|^2\right)^{\frac{p}{2}},
\end{align*}
and using the concavity of the function $s^{\frac{p}{2}}$, since $1<p<2$, we have
$$\left( t^2|\varphi_{t,L,R}'|^2+(1-t^2)|D\psi_{t,R}|^2\right)^{\frac{p}{2}}\geq t^2|\varphi_{t,L,R}'|^p+(1-t^2)|D\psi_{t,R}|^p,$$
therefore
$$-\mathcal{M}_{\Lambda,\lambda}^{+}(D^2\overline{u})+|D\overline{u}|^p\geq 0.$$
By (\ref{cotasu}), and since the function $\psi_{t,R}$ blows up at the boundary, we have that $u(\tilde{x},y)\leq \overline{u}(\tilde{x},y)$ at the boundary of the cylinder, and using comparison as in section V.1. in \cite{comparison}, we have that
$$u(\tilde{x},y)\leq t\varphi_{t,L,R}(y)+\sqrt{1-t^2}\psi_{t,R}(\tilde{x}), \quad \forall (\tilde{x},y)\in B_R^{n-1}\times (0,L).$$
Therefore, at the origin we have $u(0,y)\leq t\varphi_{t,L,R}(y)$, and translating the origin in the $\tilde{x}-$axis we have that
$$u(\tilde{x},y)\leq t\varphi_{t,L,R}(y),\quad \forall (\tilde{x},y)\in \mathbb{R}_+^n.$$
Letting $R\rightarrow \infty$, by (\ref{condicionr}) we see that $C_R$ tends to 0, therefore
\begin{equation} \label{cotacont}
u(\tilde{x},y)\leq t\varphi_{t,L}(y),\quad \forall (\tilde{x},y)\in \mathbb{R}_+^n.
\end{equation}
Now, letting $L$ goes to infinity, the function $\varphi_{t,L}$ converges to the solution of the problem
\begin{equation*}
\begin{cases}
-\Lambda \varphi_t''+t|\varphi_t'|^p=0, \mbox{ in } (0,\infty),\\
\varphi_t(0)=\dfrac{M}{t}, \quad \displaystyle\lim_{y\rightarrow \infty} \varphi_t(y)=\dfrac{l}{t}.
\end{cases}
\end{equation*}
Using (\ref{cotacont}), and $L$ tending to infinity, we have that $u(\tilde{x},y)\leq t\varphi_t(y)$ for any $t\in (0,1)$, and letting $t$ tends to 1, the function $\varphi_t$ converges to the function $v_l(y)$ defined in (\ref{onedimentional}). Then we have the inequality
\begin{equation} \label{primeradesigualdad}
u(\tilde{x},y)\leq v_l(y).
\end{equation}

\bigskip

\noindent Now we want to prove the inequality $u(\tilde{x},y)\geq v_l(y)$. Let $a\geq 0$, consider the functions $v=v_{a,R,S}(\rho)$, radial solutions of
\begin{equation} \label{radial}
\begin{cases}
-\mathcal{M}_{\Lambda,\lambda}^{+}(D^2v)+|Dv|^p=0, \mbox{ in } B_{R+S}(0)\setminus B_R(0),\\
v(R)=M, \quad v(R+S)=a,
\end{cases}
\end{equation}
and for $x=(\tilde{x},y)$ fixed with $y\in (0,S)$, the sequence $\left\lbrace v_{a,R,S}(x-x_R)\right\rbrace _{R}=\left\lbrace v_{a,R,S}(0,y+R)\right\rbrace _R$, with $x_R=(\tilde{x},-R)$. Letting $R$ tends to infinity, we can see using Lemma 3.1 in \cite{cutri} that this sequence of radial solutions converges to the unique one-dimentional solution $v_{a,S}(y)$ of the problem
\begin{equation*}
\begin{cases}
\Lambda v_{a,S}''	=|v_{a,S}'|^p \mbox{ in } (0,S),\\
v_{a,S}(0)=M, \quad v_{a,S}(S)=a.
\end{cases}
\end{equation*}
Now as $S$ goes to infinity we can see that $v_{a,S}$ converges to the function $v_a(y)$, which is the unique solution of
\begin{equation*}
\begin{cases}
\Lambda v_a''=|v_a'|^p \mbox{ in } (0,\infty),\\
v_a(0)=M,\quad \displaystyle\lim_{y\rightarrow \infty} v_a(y)=a.
\end{cases}
\end{equation*}
Since $u(\tilde{x},\cdot)\geq 0$ for all $\tilde{x}\in \mathbb{R}^{n-1}$ and $\displaystyle\lim_{y\rightarrow \infty}v_0(y)=0$, using comparison again as in \cite{comparison}, we have that $u\geq v_0(y)$, and by (\ref{cotasu}):
$$u(\tilde{x},y)\geq a_1:= \min \left[ \max \left\lbrace v_0(y),l-Cy^{-\frac{2-p}{p-1}}\right\rbrace \right],$$
Now, if we know that $u(\tilde{x},y)\geq a$ for all $(\tilde{x},y)\in \mathbb{R}_+^n$, as above, using comparison we deduce that $u(\tilde{x},y)\geq v_{a,R,S}(x- x_R)$, and letting $R$ and $S$ tend to infinity, we have that $u(\tilde{x},y)\geq v_a(y)$. Applying this to $a_1$ we deduce that $u(\tilde{x},y)\geq v_{a_1}(y)$ and then
$$u(\tilde{x},y)\geq a_2:= \min \left[ \max \left\lbrace v_{a_1}(y),l-Cy^{-\frac{2-p}{p-1}}\right\rbrace \right].$$
Then, iterating the process, we get the positive real sequence $\left\lbrace a_n\right\rbrace $ and the functions sequence $\left\lbrace v_{a_n}(y)\right\rbrace $ such that
$$u\geq v_{a_n}(y),\quad \mbox{with } a_n=\min \left[ \max \left\lbrace v_{a_{n-1}}(y),l-Cy^{-\frac{2-p}{p-1}}\right\rbrace \right] .$$
Now we letting $n$ goes to infinity, and we have that $a_n \rightarrow l$ and $v_{a_n}$ converges to $v_{l}(y)$, so we have the inequality 
$$u(\tilde{x},y)\geq v_l(y).$$
With this and (\ref{primeradesigualdad}) we conclude that $u(\tilde{x},y)=v_l(y)$.\\\\
\noindent For the case $p=2$, with a simple ODE analysis we have that the only nonnegative solution of the equation $\Lambda  v''=|v'|^2$ is $v=M$. We also note that, for $p\in (1,2)$, the solution of equation (\ref{onedimentional}) (that is $u$) is given by
$$u(\tilde{x},y)=l+\displaystyle\int_y^{\infty} \left[ \left(\dfrac{p-1}{\Lambda }\right)( t+C_{M,l}) \right]^{-\frac{1}{p-1}}dt,$$
where $C_{M,l}$ is determined by the relation
$$\displaystyle\int_0^{\infty} \left[ \left( \dfrac{p-1}{\Lambda }\right)( t+C_{M,l})\right]^{-\frac{1}{p-1}}dt=M-l.$$
\begin{flushright}
$\square$
\end{flushright}

\bigskip

\noindent \textit{Proof of Theorem \ref{teoremapmayora2}:} As in Theorem \ref{teoremapmenora2}, we have that, using Theorem \ref{teoremabernstein} applied to the function $-u$ in order to have (\ref{bernstein}), $u$ is a classic solution of equation (\ref{ecuacioneliptica2}).  Let now $x=(\tilde{x},y) \in \mathbb{R}^{n-1}\times [0,\infty)$ and $h\in \mathbb{R}^{n-1}\setminus \left\lbrace 0\right\rbrace $ fixed. We define the function $$v(\tilde{x},y)=u(\tilde{x}+h,y)-u(\tilde{x},y),\quad (\tilde{x},y)\in \mathbb{R}^{n-1}\times [0,\infty).$$
It suffices to prove that $v \equiv 0$. We assume for contradiction that
\begin{equation} \label{sigmadef}
\sigma :=\displaystyle\sup_{\mathbb{R}_+^n} v>0
\end{equation}
(the case with $\inf$ is analogous). 
Using now that $u$ is a classic solution, we have the same Bernstein estimation as in (\ref{bernsteintodopunto}), therefore $$|v(\tilde{x},y)|\leq C(n,p)|h|y^{-\frac{1}{p-1}},\quad (\tilde{x},y) \in \mathbb{R}^{n-1} \times (0,\infty).$$
It follows that $|v|\leq \frac{\sigma}{2}$ for $y\geq A$ large. Then
\begin{equation} \label{supvinA}
\sigma=\displaystyle\sup_{\mathbb{R}^{n-1}\times (0,A)} v.
\end{equation}
On the other hand, we can write $u_h=u(\tilde{x}+h,y)$ and if we define $$F(s)=|sDu_h+(1-s)Du|^p,$$ with $s\in \mathbb{R}$, we have $$F'(s)=p|sDu_h+(1-s)Du|^{p-2}(sDu_h+(1-s)Du)\cdot Dv.$$ Writing $G(\xi)=p|\xi|^{p-2}\xi$, it follows that $$|Du_h|^p-|Du|^p=F(1)-F(0)=\displaystyle\int_0^1 F'(s)ds,$$ and then $$-\mathcal{M}_{\lambda,\Lambda}^{\pm}(D^2u_h)+\mathcal{M}_{\lambda,\Lambda}^{\pm}(D^2u)=a(\tilde{x},y)\cdot Dv,$$ with $a(\tilde{x},y)=\displaystyle\int_0^1 G(sDu_h+(1-s)Du)ds$. Using now (\ref{elliptic}) we have that $$-\mathcal{M}_{\lambda,\Lambda}^{\pm}(D^2u_h)+\mathcal{M}_{\lambda,\Lambda}^{\pm}(D^2u)\geq -\mathcal{M}_{\lambda,\Lambda}^{\pm}(D^2v),$$ and then $v$ satisfies
\begin{equation} \label{ecuacionv}
-\mathcal{M}_{\lambda,\Lambda}^{\pm}(D^2v)\leq a(\tilde{x},y)\cdot Dv, \quad \mbox{ in } \mathbb{R}_+^n.
\end{equation}
By Theorem \ref{teoremabernstein}, the function $a$ is bounded in compact subsets of $\mathbb{R}_+^n$, recalling that now $u$ is a classic solution of (\ref{ecuacioneliptica2}). Therefore, by the strong maximum principle, the solution $v$ of (\ref{ecuacionv}) cannot achieve any local maximum in $\mathbb{R}_+^n$, and the function $v$ is not constant either by the fact that $v(\tilde{x},0)=0$ and (\ref{sigmadef}). Therefore, $\sigma$ in (\ref{supvinA}) is not attained and there exists a sequence $(\tilde{x}_j,y_j)\in \mathbb{R}^{n-1}\times (0,A)$ with $|x_j| \rightarrow \infty$ such that $$v(\tilde{x}_j,y_j)\rightarrow \sigma.$$
We define the sequence 
$$u_j(\tilde{x},y)=u(\tilde{x}_j+\tilde{x},y),\quad (\tilde{x},y) \in \mathbb{R}^{n-1}\times [0,\infty),$$ 
and it follows that
\begin{equation} \label{2.5} 
\displaystyle\sup_{(\tilde{x},y)\in \mathbb{R}_+^n}(u_j(\tilde{x}+h,y)-u_j(\tilde{x},y))=\displaystyle\sup_{(\tilde{x},y)\in \mathbb{R}_+^n}(u(\tilde{x}_j+\tilde{x}+h,y)-u(\tilde{x}_j+\tilde{x},y))=\displaystyle\sup_{\mathbb{R}_+^n} v=\sigma,
\end{equation}
and
\begin{equation} \label{2.6}
u_j(h,y_j)-u_j(0,y_j)=v(\tilde{x}_j,y_j) \rightarrow \sigma.
\end{equation}
Since $u_j$ is a solution of (\ref{ecuacioneliptica}) $\forall j$, by Theorem \ref{teoremabernstein} we have
\begin{equation} \label{bernsteinuj}
|Du_j(\tilde{x},y)|\leq C(n,p)y^{-\frac{1}{p-1}}, \quad \forall (\tilde{x},y)\in \mathbb{R}^{n-1}\times (0,\infty),
\end{equation}
and since $u_j(\tilde{x},0)=0$, integrating in the $y$ direction, it follows that for all $j$:
\begin{equation} \label{cotauj}
|u_j(\tilde{x},y)|\leq C(n,p)y^{\frac{p-2}{p-1}}, \quad \forall (\tilde{x},y) \in \mathbb{R}^{n-1}\times [0,\infty).
\end{equation}
As before, by Theorem 1.2 in \cite{cabre}, recalling that the functions $u_j$ are classic solutions of (\ref{ecuacioneliptica2}), we have together with (\ref{bernsteinuj}) that $u_j \in C_{loc}^{2,\alpha}(\mathbb{R}_+^n)$ for some $\alpha \in (0,1)$, then as a consequence of Arzela-Ascoli theorem we have that $\left\lbrace u_j\right\rbrace $ has a subsequence that converges uniformly in compacts sets to a classic solution $U$ of (\ref{ecuacioneliptica2}). Also, as a consequence of (\ref{cotauj}) we have $U\in C(\overline{\mathbb{R}_+^n})$ and $U(\tilde{x},0)=0$, and by (\ref{supvinA}) we also have $y_j \rightarrow y_{\infty} \in [0,A)$, then, by (\ref{2.6}) we get
\begin{equation} \label{2.8}
U(h,y_{\infty})-U(0,y_{\infty})=\sigma,
\end{equation}
which implies $y_{\infty}>0$.\\
Let now $$V(\tilde{x},y)=U(\tilde{x}+h,y)-U(\tilde{x},y),\quad (\tilde{x},y)\in \mathbb{R}^{n-1}\times [0,\infty).$$
By (\ref{2.5}) and (\ref{2.8}) we have $\sigma=\displaystyle\sup_{\mathbb{R}_+^n} V=V(0,y_{\infty})$, but $V$ satisfies
$$-\mathcal{M}_{\lambda,\Lambda}^{\pm}(D^2V)\leq A(\tilde{x},y)\cdot DV,$$
where $A(\tilde{x},y)=\displaystyle\int_0^1 G(sDU(\tilde{x}+h,y)+(1-s)DU(\tilde{x},y))ds$ is bounded on compact subsets of $\mathbb{R}_+^n$. Finally, this contradicts the strong maximum principle.
\begin{flushright}
$\square$
\end{flushright}
\bigskip
As we present in section 1, as an application to their result, Filipucci, Pucci and Souplet in \cite{fps} proved an improvement of the Bernstein estimation in \cite{lions}. The result establishes a more precise Bernstein constant when the solution $u$ is zero at the border of $\Omega$, which is also valid for an $L^q$-viscosity solution of the Pucci equation (\ref{ecuacionelipticanohomogenea}) in the case $p>2$. As in the beginning of the proof of Theorem \ref{teoremapmenora2}, considering $\tilde{f}(x)=|Du|^p+f(x)$, it can be proved that the solution $u$ is classic. \\

\noindent The Theorem \ref{teoremabernsteinmejorado} is an application to the Bernstein estimation in Theorem \ref{teoremabernstein} and the Liouville-type result from Theorem \ref{teoremapmayora2}. In the proof we follow the ideas from \cite{fps}.\\

\noindent \textit{Proof of the Theorem \ref{teoremabernsteinmejorado}:} We will use, as mentioned above, that $u$ is a classic solution. Let $c_p=\left(\dfrac{\lambda}{p-1}\right)^{\frac{1}{p-1}}$, and assume for contradiction that there exist $c>c_p$ and sequences $\left\lbrace f_j\right\rbrace ,\left\lbrace u_j\right\rbrace $ and $\left\lbrace x_j\right\rbrace $, $f_j$ Lipschitz in $\Omega$ with $\|f_j\|_{W^{1,\infty}}\leq M$, such that
\begin{equation} \label{absurdo}
\begin{cases}
-\mathcal{M}_{\lambda,\Lambda}^{\pm}(D^2u_j)=|Du_j|^p+f_j(x), \quad \mbox{in } \Omega,\\
d_{\partial \Omega}(x_j)\rightarrow 0,\\
d_{\partial \Omega}^{\frac{1}{p-1}}(x_j)|Du_j(x_j)|\geq c.
\end{cases}
\end{equation}
Let $P(x)$ the projection of $x$ onto $\partial \Omega$, and
$$z_j:=P(x_j),\quad \alpha_j:=d_{\partial \Omega}(x_j)=|x_j-z_j|.$$
By extracting a subsequence if necessary, we assume that $z_j \rightarrow a \in \partial \Omega$, and without loss of generality we may assume that $a=0$, therefore $x_j \rightarrow 0$ and we have that
\begin{equation} \label{vectornormal}
\nu_j:=\nu(z_j) \rightarrow e_n.
\end{equation}
Let now $v_j(y)=\alpha_j^{\frac{2-p}{p-1}}u_j(z_j+\alpha_j y)$, we have that
\begin{align*}
Dv_j(y)&=\alpha_j^{\frac{1}{p-1}}Du_j(z_j+\alpha_j y),\\
\mathcal{M}_{\lambda,\Lambda}^{\pm}(D^2v_j(y))&=\alpha_j^{\frac{p}{p-1}}\mathcal{M}_{\lambda,\Lambda}^{\pm}(D^2u_j(z_j+\alpha_j y)).
\end{align*}
If we define $\Omega_j:=\alpha_j^{-1}(\Omega-z_j)$, we have that the function $v_j$ satisfies the equation
$$\mathcal{M}_{\lambda,\Lambda}^{\pm}(D^2v_j)+|Dv_j|^p=\alpha_j^{\frac{p}{p-1}}[\mathcal{M}_{\lambda,\Lambda}^{\pm}(D^2u_j)+|Du_j|^p](z_j+\alpha_j y),\quad \mbox{in } \Omega_j,$$
therefore
$$-\mathcal{M}_{\lambda,\Lambda}^{\pm}(D^2v_j)=|Dv_j|^p+\tilde{f}_j(y), \quad \mbox{in } \Omega_j,$$
where $\tilde{f}_j(y)=\alpha_j^{\frac{p}{p-1}}f_j(z_j+\alpha_j y).$ We note that $\Omega_j$ converges to the half-space $\mathbb{R}_+^n$ as $j\rightarrow \infty$.\\
\noindent Now, using Bernstein estimation from Theorem \ref{teoremabernstein} it follows that
$$|Du_j(x)|\leq Cd_{\partial \Omega}^{-\frac{1}{p-1}}(x) \quad \mbox{in } \Omega,$$
and
$$|u_j(x)|\leq Cd_{\partial \Omega}^{\frac{p-2}{p-1}}(x) \quad \mbox{in } \Omega,$$
with $C=C(n,p,M)>0$, independent of $j$. Now we need to find a uniform bound for $v_j$ and $Dv_j$. Letting $$Q_{R,\epsilon}=B_R \cap \left\lbrace y_n>\epsilon\right\rbrace,$$ as in the proof of Proposition 3.1 in \cite{fps}, it can be proved that $Q_{R,\epsilon} \subset \Omega_j$ for all $j\geq j_0$, where $j_0$ depends on $R$ and $\epsilon$, and
\begin{equation} \label{cotadistanciay}
\dfrac{1}{2}y_n\leq d_{\partial \Omega_j}(y)\leq 3y_n, \quad \forall y\in Q_{R,\epsilon}, \forall j\geq j_0.
\end{equation}
Using now the Bernstein estimation from Theorem \ref{teoremabernstein}, we have that, for all $j\geq j_0$ and $y\in Q_{R,\epsilon}$,
\begin{equation*}
|Dv_j(y)|=\alpha_j^{\frac{1}{p-1}}|Du_j(z_j+\alpha_j y)|\leq C\alpha_j^{\frac{1}{p-1}}d_{\partial \Omega}^{-\frac{1}{p-1}}(z_j+\alpha_j y),
\end{equation*}
and since $d_{\partial \Omega}(z_j+\alpha_j y)=\alpha_j d_{\partial \Omega_j}(y)$, using (\ref{cotadistanciay}), we have
\begin{equation} \label{cotadvjenq}
|Dv_j(y)|\leq C2^{\frac{1}{p-1}}y_n^{-\frac{1}{p-1}}.
\end{equation}
Similarly, we have the estimation
\begin{equation} \label{cotavjenq}
|v_j(y)|\leq C3^{\frac{p-2}{p-1}}y_n^{\frac{p-2}{p-1}}.
\end{equation}
Using now interior elliptic $L^q$ estimates, we have that the sequence $\left\lbrace  v_j\right\rbrace _j$ is precompact in $W^{2,q}(Q_{R,\epsilon})$, and using a diagonal procedure, we deduce that some subsequence of $(v_j)_j$ converges in $W^{2,q}(Q_{R,\epsilon})$ for each $q,R,\epsilon>0$, to a strong solution $V(y)\in C^2(\mathbb{R}_+^n)\cap C(\overline{\mathbb{R}_+^n})$ of
\begin{equation*}
\begin{cases}
-\mathcal{M}_{\lambda,\Lambda}^{\pm}(D^2V)=|DV|^p, \quad y\in \mathbb{R}_+^n,\\
V(y)=0, \quad y\in \partial \mathbb{R}_+^n.
\end{cases}
\end{equation*}
It follows from Theorem \ref{teoremapmayora2} that either $V=0$ or $V(y)=U_{\hat{c}}(y_n)$ for some $\hat{c}\geq 0$, and by (\ref{solucionfundamentalpmayora2}) we have
$$|DV(y)|=c_p(y_n+\hat{c})^{-\frac{1}{p-1}}\leq c_py_n^{-\frac{1}{p-1}}.$$
Finally, with this fact and using (\ref{vectornormal}):
\begin{equation*}
\displaystyle\lim_{j\rightarrow \infty}d_{\partial \Omega}^{\frac{1}{p-1}}|Du_j(x_j)|=\displaystyle\lim_{j\rightarrow \infty}\,\vline Dv_j\left(\frac{x_j-z_j}{\alpha_j}\right)\, \vline\,=\displaystyle\lim_{j\rightarrow \infty}|Dv_j(\nu_j)|=|DV(e_n)|\leq c_p,
\end{equation*}
which is a contradiction with (\ref{absurdo}), being $c>c_p$.
\begin{flushright}
$\square$
\end{flushright}

\section{Appendix}

In this section we prove Lipschitz regularity for any continuous $L^q$-viscosity solution of the equation (\ref{ecuacion1}), and in addition, we assume that there exist positive constants $\gamma_0$ and $\gamma_1$ such that, for all $\xi, \eta \in \mathbb{R}^n$,
\begin{equation} \label{3.3}
\begin{array}{cc}
&g(\xi)\geq \gamma_0 |\xi|^p,\quad p>1,\\
&|g(\xi)-g(\eta)|\leq \gamma_1 (|\xi|^{p-1}+|\eta|^{p-1})|\xi-\eta|.
\end{array}
\end{equation}

\bigskip

\begin{teo} \label{teoremabernstein}
Let $\Omega \subset \mathbb{R}^n$ open and bounded, with Lipschitz boundary and satisfy the uniform interior sphere condition. Let $p>1$ and $u\in C(\Omega)$ a bounded $L^q$-viscosity solution of (\ref{ecuacion1}) in $\Omega$, with $n\leq q$. Assume that (\ref{3.3}) hold true, and that the function $f$ is Lipschitz and bounded in $\Omega$, and $d_{\partial \Omega}^{\frac{p}{p-1}}f$ is  also bounded in $\Omega$, where $d_{\partial \Omega}(x)=dist(x,\partial \Omega)$. Then $u$ is locally Lipschitz continuous and there exists a constant $C>0$, depending on $p$ and $\|d_{\partial \Omega}^{\frac{p}{p-1}}f\|_{L^{\infty}(\Omega)}$, such that $$|Du(x)|\leq \dfrac{C}{d_{\partial \Omega}(x)^{\frac{1}{p-1}}} \mbox{   for } a.e.\,\,\, x\in \Omega.$$
\end{teo}

For the proof of Theorem \ref{teoremabernstein}, we present a simplification of the work done in \cite{bernstein}, following the ideas from \cite{clp}. The work of Capuzzo Dolcetta, Leoni and Porreta in \cite{clp} (specifically, section 3) is for the equation
$$-tr(A(x)D^2u)+\alpha u+H(x,Du)=0, \quad \mbox{in  } \Omega,$$
where $A$ is a bounded and Lipschitz continuous map, and $\alpha \geq 0$. As in \cite{clp}, the assumption of non-degeneracy is not required for Theorem \ref{teoremabernstein}, considering that the proof is valid in the case $\lambda=0$. The technique consists in a variable bending, using a $C^2$- test function depending on the distance to the boundary, with the calculation of the Pucci operators from derivatives of the test function,  which are non-linear operators.\\\\
\noindent In the first case we are going to prove the result in $B(0,1)$ with exponent $\frac{p}{p-1}$ instead of $\frac{1}{p-1}$ as in the theorem. Then, when rescaling to a domain $\Omega$, we will improve the Bernstein bound.\\

\noindent \textit{Proof of Theorem \ref{teoremabernstein}:} \textbf{First case:} $\Omega=B(0,1)=B$. We choose a smooth monotone radial function $d\in C^2(\overline{B})$ satisfying the properties: 
$$\left\lbrace \begin{array}{c}
d(x)=d_{\partial B}(x) \mbox{ if } d_{\partial B}\leq \dfrac{1}{2},\\
\dfrac{d_{\partial B}(x)}{2}\leq d(x)\leq d_{\partial B}(x) \,\,\, \forall x\in \overline{B},\\
|Dd(x)|\leq 1, D^2d(x)\leq 0 \,\,\, \forall x\in \overline{B}.
\end{array} \right. $$
Now consider the function $$\Phi (x,y)=k|x-y|\varphi (x,y),$$
where 
\begin{equation} \label{3.5}
\varphi (x,y)=\dfrac{1}{d(y)^{\gamma}}\left[ L+\left( \dfrac{|x-y|}{d(x)}\right) ^{\beta}\right] .
\end{equation}
The constants $\gamma, \beta$ y $L$ are positive and will be fixed later. We want to prove that, choosing $k>>1$ we have $$u(x)-u(y)-\Phi (x,y)\leq 0 \,\,\, \forall x,y \in B.$$
Suppose by contradiction, that $w(x,y)=u(x)-u(y)-\Phi(x,y)>0$ at some point $(x,y)\in B\times B$, this says that $w$ has a maximum point at $(x,y) \in B\times B$, with $x\neq y$. Observing that the function $\Phi$ is $C^2(B)$, from \cite{users}, specifically \textit{Theorem 3.2}, we have that for every $\epsilon >0$ there exist matrices  $X=X(\epsilon)$, $Y(\epsilon)\in \mathbb{S}_n$ such that
\begin{eqnarray} \label{3.6}
&(D_x \Phi(x,y),X) \in \overline{J}^{2,+}u(x)_B, (-D_y \Phi(x,y),Y)\in \overline{J}^{2,-}_Bu(y),\\
&\notag -\left( \dfrac{1}{\epsilon}+\| D^2 \Phi(x,y)\| \right) I_{2n}\leq \left( \begin{array}{cc}
X&0\\
0&-Y
\end{array}\right) \leq D^2 \Phi(x,y)+\epsilon (D^2 \Phi (x,y))^2,
\end{eqnarray}
where $\overline{J}^{2,\pm}$ denotes the closure of the second order super(sub)-jet. Using the definition from sub and super-jet in \cite{users} we have:
\begin{eqnarray} \label{3.7}
&-\mathcal{M}_{\lambda,\Lambda}^{\pm}(X) +g(D_x \Phi)\leq f(x),\\
&-\mathcal{M}_{\lambda,\Lambda}^{\pm}(Y)+g(-D_y\Phi) \geq f(y),
\end{eqnarray}
and multiplying the first equation by (1+t), for fixed $t>0$: $$\mathcal{M}_{\lambda,\Lambda}^{\pm}(Y)-(1+t)\mathcal{M}_{\lambda,\Lambda}^{\pm}(X) +(1+t)g(D_x\Phi)-g(-D_y\Phi)\leq (1+t)f(x)-f(y).$$
Using (\ref{elliptic}) and the fact that $\mathcal{M}_{\lambda,\Lambda}^+(M)=-\mathcal{M}_{\lambda,\Lambda}^-(-M)$ we have
\begin{equation} \label{3.8}
-\mathcal{M}_{\lambda,\Lambda}^+((1+t)X-Y)+(1+t)g(D_x\Phi)-g(-D_y\Phi)\leq (1+t)f(x)-f(y),
\end{equation}
and using (\ref{3.3}):
\begin{equation} \label{3.9}
\gamma_0 t|D_x\Phi|^p\leq \mathcal{M}_{\lambda,\Lambda}^+((1+t)X-Y)+g(-D_y\Phi)-g(D_x\Phi)+(1+t)f(x)-f(y).
\end{equation}
In the calculation of derivatives of $\Phi$, we have: 
\begin{equation} \label{3.12}
D_x\Phi=\dfrac{k}{d(y)^{\gamma}}\left[ \left( L+(1+\beta)\left( \dfrac{|x-y|}{d(x)}\right) ^{\beta}\right) \dfrac{(x-y)}{|x-y|}-\beta \left( \dfrac{|x-y|}{d(x)}\right) ^{\beta +1}Dd(x)\right].
\end{equation}
Denoting $\xi=\dfrac{|x-y|}{d(x)}$, and using the notation $\widehat{v}=\frac{v}{|v|}$ we have the estimation 
$$|D_x\Phi|^2\geq \dfrac{k^2}{d(y)^{2\gamma}}\left[ \dfrac{1}{2}(L+(1+\beta)\xi^{\beta})^2+\beta^2 |Dd(x)|^2 \xi^{2(\beta+1)}-2\beta^2\xi^{2\beta}\right] .$$
Since $|Dd(x)|=1$ if $d(x)\leq \frac{1}{2}$ and $\xi=\frac{|x-y|}{d(x)}\leq 4$ if $d(x)\geq \frac{1}{2}$, in both cases there exists $c>0$, depending only on $\beta$, such that:
$$|D_x\Phi|^2\geq \dfrac{k^2}{d(y)^{2\gamma}}\left[ \dfrac{1}{2}(L+(1+\beta)\xi^{\beta})^2+\dfrac{\beta^2}{2}\xi^{2(\beta+1)}-c\right] .$$
We choose $L$ sufficiently large such that $$|D_x\Phi|^2\geq c\dfrac{k^2}{d(y)^{2\gamma}}\left[ (L+(1+\beta)\xi^{\beta})^2+\beta^2\xi^{2(\beta+1)}\right], $$ and therefore
\begin{equation} \label{3.14}
|D_x\Phi|\geq c\dfrac{k}{d(y)^{\gamma}}(L+\xi^{\beta})(1+\xi)=ck\varphi (1+\xi).
\end{equation}
Now, let $A\in \mathcal{A}_{\lambda,\Lambda}$ and the non negative matrix 
\begin{equation} \label{3.10}
\mathbf{A_t}=\left( \begin{array}{cc}
(1+t)A&\sqrt{1+t}A\\
\sqrt{1+t}A&A
\end{array}\right) .
\end{equation}
Multiplying the right inequality in (\ref{3.6}) by $\mathbf{A_t}$, taking traces and then taking supreme over $A\in \mathcal{A}_{\lambda,\Lambda}$:
$$\mathcal{M}_{\lambda,\Lambda}^+((1+t)X-Y) \leq \displaystyle\sup_{A\in \mathcal{A}_{\lambda,\Lambda}} tr(\mathbf{A}_tD^2\Phi)+\epsilon \displaystyle\sup_{A\in \mathcal{A}_{\lambda,\Lambda}}tr(\mathbf{A_t}(D^2\Phi)^2) ,$$
and from (\ref{3.9}) we have: 
\begin{equation} 
t|D_x\Phi|^p\leq \displaystyle\sup_{A\in \mathcal{A}_{\lambda,\Lambda}} tr(\mathbf{A}_tD^2\Phi)+\epsilon \displaystyle\sup_{A\in \mathcal{A}_{\lambda,\Lambda}}tr(\mathbf{A_t}(D^2\Phi)^2) +|D_y\Phi|^p-|D_x\Phi|^p.
\end{equation}
Now letting $\epsilon$ tend to 0:
\begin{equation}\label{3.11}
t|D_x\Phi|^p\leq \displaystyle\sup_{A\in \mathcal{A}_{\lambda,\Lambda}} tr(\mathbf{A}_tD^2\Phi) +|D_y\Phi|^p-|D_x\Phi|^p.
\end{equation}
In calculating $D^2\Phi$, multiplying by $\mathbf{A}_t$ and taking traces we get the expression
\begin{align*} 
tr(\mathbf{A_t}D^2\Phi(x,y))&=\dfrac{k}{d(y)^{\gamma}}\left\lbrace \dfrac{(L+(1+\beta)\xi^{\beta})}{|x-y|}(2-2\sqrt{1+t}+t)\mathcal{M}_{\lambda,\Lambda}^+(B)\right. \\
&+\dfrac{\beta (1+\beta)\xi^{\beta}}{|x-y|}(2-2\sqrt{1+t}+t)\mathcal{M}_{\lambda,\Lambda}^+(T)\\
&+\dfrac{\beta (1+\beta)\xi^{\beta}}{d(x)}\sqrt{1+t}\,\mathcal{M}_{\lambda,\Lambda}^+(-\sqrt{1+t}C-(\sqrt{1+t}-2)C^T)\\
&+\dfrac{\gamma (L+(1+\beta)\xi^{\beta})}{d(y)}(\sqrt{1+t}-1)\mathcal{M}_{\lambda,\Lambda}^+(-D-D^T)\\
&+\dfrac{\beta (1+\beta)\xi^{\beta+1}}{d(x)}(1+t) \mathcal{M}_{\lambda,\Lambda}^+(Dd(x)\otimes Dd(x))\\
&+\dfrac{\gamma (\gamma+1)|x-y|}{d(y)^2}(L+\xi^{\beta})\mathcal{M}_{\lambda,\Lambda}^+(Dd(y)\otimes Dd(y))\\
&+\dfrac{\beta \gamma \xi^{\beta+1}}{d(y)}\sqrt{1+t}\left[ \mathcal{M}_{\lambda,\Lambda}^+(Dd(x)\otimes Dd(y))+\mathcal{M}_{\lambda,\Lambda}^+(Dd(y)\otimes Dd(x))\right] \\
&+\left. \beta \xi^{\beta+1}(1+t)\mathcal{M}_{\lambda,\Lambda}^+(-D^2d(x))+\dfrac{\gamma |x-y|}{d(y)}(L+\xi^{\beta}) \mathcal{M}_{\lambda,\Lambda}^+(-D^2d(y))\right\rbrace .
\end{align*}
where $B,T,C$ and $D$ are the matrices defined by $B=I-\widehat{x-y}\otimes \widehat{x-y}$, $T=\widehat{x-y}\otimes \widehat{x-y}$, $C=\widehat{x-y}\otimes Dd(x)$ and $D=\widehat{x-y}\otimes Dd(y)$, with $\otimes$ the Kronecker product between two vectors. Using the fact that for any vectors $v_1,v_2 \in \mathbb{R}^n$ the only non trivial eigenvalue from the matrix $v_1 \otimes v_2$ is $v_1\cdot v_2$ and property in (\ref{elliptic}), the next estimation is valid:
$$\mathcal{M}_{\lambda,\Lambda}^+(-\sqrt{1+t}C-(\sqrt{1+t}-2)C^T)=-2(\sqrt{1+t}-1)\mathcal{M}_{\lambda,\Lambda}^-(C),$$
and since $\mathcal{M}_{\lambda,\Lambda}^-(C)=\mathcal{M}_{\lambda,\Lambda}^-(\widehat{x-y}\otimes Dd(x))=C_{\lambda, \Lambda}(\widehat{x-y})\cdot Dd(x)$, we have $$\mathcal{M}_{\lambda,\Lambda}^+(-\sqrt{1+t}C-(\sqrt{1+t}-1)C^T)\leq -2C_{\lambda, \Lambda}(\sqrt{1+t}-1)(\widehat{x-y})\cdot Dd(x).$$
On the other hand:
$$\mathcal{M}_{\lambda,\Lambda}^+(-D-D^T) \leq -\left[ \mathcal{M}_{\lambda,\Lambda}^-(D)+\mathcal{M}_{\lambda,\Lambda}^-(D^T)\right] =-2\mathcal{M}_{\lambda,\Lambda}^-(D),$$ and since $\mathcal{M}_{\lambda,\Lambda}^-(D)=\mathcal{M}_{\lambda,\Lambda}^-(\widehat{x-y}\otimes Dd(y))=C_{\lambda, \Lambda} (\widehat{x-y})\cdot Dd(y)$ we have $$\mathcal{M}_{\lambda,\Lambda}^+(-D-D^T)\leq -2C_{\lambda, \Lambda}(\widehat{x-y})\cdot Dd(y).$$
\bigskip
Now observe that for all $t>0$ we have 
\begin{equation*}
\begin{cases}
(2-2\sqrt{1+t}+t)\leq t^2 \leq t^2+|x-y|^2,\\
2\sqrt{1+t}(\sqrt{1+t}-1)\leq 2t \leq 2t+|x-y|,
\end{cases}
\end{equation*}
and we obtain then, 
\begin{eqnarray}  \label{3.16}
\quad \operatorname{tr}(\mathbf{A_t}D^2\Phi) &\leq \dfrac{ck}{d(y)^{\gamma}}\left\lbrace \dfrac{(L+\xi^{\beta})}{|x-y|}(t^2+|x-y|^2)+\left( \dfrac{\xi^{\beta}}{d(x)}+\dfrac{1+\xi^{\beta}}{d(y)}\right) (t+|x-y|)\right. \\
&\notag \left. +\xi^{\beta+1}\left( \dfrac{1}{d(x)}+\dfrac{1}{d(y)}\right) t+\left( \dfrac{\xi^{\beta+1}}{d(x)}+\dfrac{\xi^{\beta+1}}{d(y)}+\dfrac{(1+\xi^{\beta})|x-y|}{d(y)^2}\right) \right\rbrace .
\end{eqnarray}
We also have the estimation:
$$\dfrac{\xi^{\beta}}{d(x)}\leq 2\dfrac{(1+\xi^{\beta+1})}{d(y)},$$
and with this, we obtain
$$\dfrac{\xi^{\beta+1}}{d(x)}+\dfrac{\xi^{\beta+1}}{d(y)}+\dfrac{(1+\xi^{\beta})|x-y|}{d(y)^2}\leq c|x-y|\dfrac{(1+\xi^{\beta+2})}{d(y)^2},$$
and
$$\left( \dfrac{\xi^{\beta}}{d(x)}+\dfrac{1+\xi^{\beta}}{d(y)}\right) (t+|x-y|)+\xi^{\beta+1}\left( \dfrac{1}{d(x)}+\dfrac{1}{d(y)}\right) t\leq c\dfrac{(1+\xi^{\beta+1})}{d(y)}|x-y|+ct\dfrac{(1+\xi^{\beta+2})}{d(y)}.$$
Using this in (\ref{3.16}) and by (\ref{3.11}) we have:
\begin{align*}
\gamma_0 t|D_x\Phi|^p &\leq g(-D_y\Phi)-g(D_x\Phi)+(1+t)f(x)-f(y)+c\dfrac{k}{d(y)^{\gamma}}\left[ \dfrac{(L+\xi^{\beta})}{|x-y|}(t^2+|x-y|^2)\right. \\
&\left. +t\dfrac{(1+\xi^{\beta+2})}{d(y)}+\dfrac{(1+\xi^{\beta+1})}{d(y)}|x-y|+\dfrac{(1+\xi^{\beta+2})}{d(y)^2}|x-y|\right] .
\end{align*}
Considering now $\gamma \geq \dfrac{1}{p-1}$ and $\beta \geq \dfrac{2-p}{p-1}$ we have $\gamma p\geq \gamma +1$ and $(\beta+1)p\geq \beta+2$, and by estimation of $|D_x\Phi|^p$ in (\ref{3.14}):
$$|D_x\Phi|^p \geq ck^p \dfrac{(1+\xi^{\beta+1})^p}{d(y)^{\gamma p}}\geq ck^p \dfrac{(1+\xi^{\beta+2})}{d(y)^{\gamma+1}}.$$
Since $p>1$, taking $k^{p-1}>1$ we have:
$$ctk\dfrac{(1+\xi^{\beta+2})}{d(y)^{\gamma+1}}\leq ct k^p\dfrac{(1+\xi^{\beta+2})}{d(y)^{\gamma+1}}\leq \dfrac{t}{2}|D_x\Phi|^p,$$ and therefore:
\begin{align*}
\dfrac{t}{2}|D_x\Phi|^p-ckt^2\dfrac{(L+\xi^{\beta})}{|x-y|d(y)^{\gamma}}&\leq tf(x)+g(-D_y\Phi)-g(D_x\Phi)+f(x)-f(y)\\
&+ck\dfrac{|x-y|}{d(y)^{\gamma}}\left[ (1+\xi^{\beta})+\dfrac{1+\xi^{\beta+1}}{d(y)}+\dfrac{1+\xi^{\beta+2}}{d(y)^2}\right] .
\end{align*}
Now we choose optimal $t$ to maximize the left side: 
\begin{equation} \label{3.19}
t=c\dfrac{|D_x\Phi|^p}{k\frac{(L+\xi^{\beta})}{|x-y|d(y)^{\gamma}}}=c\dfrac{|D_x\Phi|^p}{k}\cdot \dfrac{|x-y|}{\varphi}.
\end{equation}
Then
\begin{align*}
\dfrac{|D_x\Phi|^{2p}|x-y|}{k\varphi}&\leq c\left\lbrace \dfrac{|D_x\Phi|^p}{k}\cdot \dfrac{|x-y|}{\varphi}\|f\|+ g(-D_y\Phi)-g(D_x\Phi)+f(x)-f(y)\right. \\
&\left. +k|x-y|\left[ \varphi+\dfrac{1+\xi^{\beta+1}}{d(y)^{\gamma+1}}+\dfrac{1+\xi^{\beta+2}}{d(y)^{\gamma+2}}\right]\right\rbrace  ,
\end{align*}
and therefore:
\begin{equation} \label{3.20}
|D_x\Phi|^{2p}\leq c\left\lbrace |D_x\Phi|^p \|f\|+ k\dfrac{\varphi}{|x-y|}\left[ g(-D_y\Phi)-g(D_x\Phi)+f(x)-f(y)\right] +k^2\varphi\left[ \varphi+\dfrac{1+\xi^{\beta+2}}{d(y)^{\gamma+2}}\right]\right\rbrace  .
\end{equation}
On the estimation of the term $g(-D_y\Phi)-g(D_x\Phi)$, using hypothesis in (\ref{3.3}) we can get
$$g(-D_y\Phi)-g(D_x\Phi)\leq ck^p \dfrac{(1+\xi^{\beta+1})^p}{d(y)^{\gamma p+p}}|x-y|.$$
Using this in (\ref{3.20}), with the fact that $f$ is Lipschitz we have:
$$|D_x\Phi|^{2p}\leq c\left\lbrace |D_x\Phi|^p\|f\|+k\varphi+ k^{p+1}\varphi \dfrac{(1+\xi^{\beta+1})^p}{d(y)^{\gamma p+p}}+k^2\varphi \left( \varphi+\dfrac{1+\xi^{\beta+2}}{d(y)^{\gamma+2}}\right) \right\rbrace ,$$
and choosing now $\gamma=\dfrac{p}{p-1}$, using (\ref{3.14}), we can find the estimation:
$$k^{p+1}\varphi \dfrac{(1+\xi^{\beta+1})^p}{d(y)^{\gamma p+p}}\leq \dfrac{c}{k^{p-1}}|D_x\Phi|^{2p}.$$
Also, since $2(\gamma+1)\leq 2\gamma p$ y $\beta+2\leq (\beta+1)p$ we have:
$$k^2\varphi \dfrac{(1+\xi^{\beta+2})}{d(y)^{\gamma+2}}\leq c\dfrac{|D_x\Phi|^{2p}}{k^{2(p-1)}}.$$
Therefore:
\begin{align*}
|D_x\Phi|^{2p}&\leq c\left\lbrace |D_x\Phi|^p\|f\|+k\varphi +k^2\varphi^2 +\dfrac{|D_x\Phi|^{2p}}{k^{p-1}}+\dfrac{|D_x\Phi|^{2p}}{k^{2(p-1)}}\right\rbrace \\
&\leq c\left\lbrace |D_x\Phi|^p+|D_x\Phi|+|D_x\Phi|^2+\dfrac{|D_x\Phi|^{2p}}{k^{p-1}}+\dfrac{|D_x\Phi|^{2p}}{k^{2(p-1)}}\right\rbrace ,
\end{align*}
which is a contradiction taking $k>>1$, by (\ref{3.14}) . Therefore there is a constant $M>0$, depending on $p$ and $f$, such that
$$u(x)-u(y)\leq M\dfrac{|x-y|}{d_{\partial B}(y)}\left[ 1+\left( \dfrac{|x-y|}{d_{\partial B}(x)}\right) ^{\beta}\right] \,\, \forall x,y\in B.$$
Changing the roles of $x$ and $y$:
$$|u(x)-u(y)|\leq M\dfrac{|x-y|}{(d_{\partial B}(y)\wedge d_{\partial B}(x))^{\gamma}}\left[ 1+\left( \dfrac{|x-y|}{d_{\partial B}(x) \wedge d_{\partial B}(y)}\right) ^{\beta}\right], $$ for all $x,y \in B$, so that $u$ is locally Lipschitz in $B$ and satisfies
$$|Du(x)|\leq \dfrac{M}{d_{\partial B}(x)^{\gamma}}, \mbox{   for } a.e. \,\,\,x\in B,$$
with $\gamma =\dfrac{p}{p-1}$.\\\\
\textbf{Second case:} In order to have a Bernstein estimation of the solution in $\Omega$, we do the following rescale. Let $x_0\in \Omega$ a differentiability point for $u$ , and let $r=\dfrac{d_{\partial \Omega}(x_0)}{2}$, we define
$$v(x)=r^{\frac{2-p}{p-1}}u(x_0+rx),\,\, x\in B(0,1).$$
For this function we have $$Dv(x)=r^{\frac{1}{p-1}}Du(x_0+rx),$$ and $$-\mathcal{M}_{\lambda,\Lambda}^+(D^2v(x))=-r^{\frac{p}{p-1}}\mathcal{M}_{\lambda,\Lambda}^+(D^2u(x_0+rx)),$$ therefore, $v$ is a viscosity solution from the equation
$$-\mathcal{M}_{\lambda,\Lambda}^{\pm}(D^2v)+g(Dv)=r^{\frac{p}{p-1}}f(x_0+rx), \quad x\in B,$$
and by choose from $r$ we have $r\leq d_{\partial \Omega}(x_0+rx)\leq 3r$, hence
$$\|r^{\frac{p}{p-1}}f\|_{L^{\infty}(B)}=\|\frac{d_{\partial \Omega}^{\frac{p}{p-1}}(x_0)}{2}f(x_0+rx)\|_{L^{\infty}(B)}\leq \| d_{\partial \Omega}^{\frac{p}{p-1}}(x_0+rx)f(x_0+rx)\|_{L^{\infty}(B)}=\|d_{\partial \Omega}^{\frac{p}{p-1}}f\|_{L^{\infty}(\Omega)},$$
therefore $$|Dv(x)|\leq \dfrac{C}{d_{\partial B}(x)^{\frac{p}{p-1}}},$$ and evaluating at $x=0$:
\begin{equation}
|Du(x_0)|\leq \dfrac{C}{r^{\frac{1}{p-1}}}=\dfrac{C}{d_{\partial \Omega}(x_0)^{\frac{1}{p-1}}},
\end{equation}
where $C$ depends on $p$ and $\| d_{\partial \Omega}^{\frac{p}{p-1}}f\|_{L^{\infty}(\Omega)}$.
\begin{flushright}
$\square$
\end{flushright}
\bigskip
\begin{rem}
The above rescale also applies to unbounded domain $\mathbb{R}_+^n$, taking $r=\frac{x_{0,n}}{2}$, where $x_{0,n}$ is the $n$-th coordinate from a point $x_0$.
\end{rem}
\bigskip
{\bf Acknowledgements} 
R. F. was supported by ANID-Magister Nacional Grant \# 22190374.
A. Q. was partially supported by FONDECYT Grant \# 1190282 and Programa Basal, CMM. U. de Chile.

\end{document}